\documentclass{article}
\usepackage[utf8]{inputenc}
\usepackage{amsmath,amssymb,amsthm}

\renewcommand\le{\leqslant}
\renewcommand\ge{\geqslant}
\newcommand\eps{\varepsilon}

\newtheorem*{theoremA}{Theorem}
\newcommand{\Z}{\mathbb{Z}}

\title{The algorithm for the recovery of integer vector via linear measurements.}
\author{K.S. Ryutin\thanks{Laboratory ``High–dimensional approximation and applications'', Department of Mechanics and Mathematics,
Moscow State University, Moscow, Russia; Email: kriutin@yahoo.com. Research supported by the grant of the Government of the Russian Federation
(project 14.W03.31.0031).}}

\begin{document}
\maketitle
\begin{abstract}
In this paper we continue the studies on the integer sparse recovery problem that was introduced in \cite{FKS} and studied in \cite{K},\cite{KS}. We provide an algorithm for the recovery of an unknown sparse integer vector for the measurement matrix described in \cite{KS} and estimate the number of arithmetical operations.
\end{abstract}

We consider the problem of recoverying some sparse integer vector via a small number of linear measurements.   It is related to the compressed sensing theory (see \cite{FR}).  The problem under consideration appeared as a continuation of studies in papers  \cite{FKS} (it introduced the integer valued compressed sensing problem) and  \cite{K}, \cite{KS}. 

 In \cite{KS} a very natural construction of a measurement  integer  matrix with a good control on the absolute values of its elements that permits to recover any sparse vector was proposed. But the question of a good algorithm for recovery was not addressed. Let  $\Phi=\Phi_{m\times p}$  be a matrix with elements $\phi_{ij}\in \Z,$ with $\phi_{ij}=k_j j^{i-1} \pmod{p}, 1\le j\le p,1\le i\le m,$  $p$ -- prime number and $m\le p,$  $k_j$ -- some fixed integers $1\le k_j<p$ (so, a representative of the residue class is chosen). In \cite{KS} the elements were chosen such that  $|\phi_{ij}|\le p^{1-1/m}$. We want to recover a vector  $x\in \Z^p,$ using the following information: the number of its non-zero coordinates $\|x\|_0=s\le m/2$ and we are given  $y=(y_0,\dots,y_{m-1})=\Phi x\in \Z^m.$ Let $I=\{j_1,\dots , ,j_s\}$  be the support of $x$, and  $x_j \in \Z \setminus\{0\}$,  $j\in I$ are the corresponding coordinates of  $x$.    In what follows we identify  the indices in $I$ with elements in the field  $F=\mathbb{F}_p.$
The  problem under consideration boils down to finding from the following system of equations
\begin{equation}
\label{sys1}
\sum_{j\in I} k_j j^{l}x_j=y_{l},\,\,0\le l \le m-1. 
\end{equation}
of the  set  $I\subset \{1,\dots, p\}$ and coefficients $x_j\in \Z.$  
If the set  $I$ is known then the  coefficients    $x_j\in\mathbb{Z}$ can be determined through the system of  $m$ linear equations in  $s$ variables, given by a matrix with all elements satisfying the estimate  $|\phi_{ij}|\le p^{1-1/m}.$  Different effective algorithms were developed for this problem.

Therefore our goal is to find the set  $I.$
 We do not know the exact value of   $s$ and the  complexity of the algorithm is measured in terms of  $m,p,M=\max{|x_j|}.$ 

We will make use of some well--known algorithms in finite fields. 
Let  $Z(t)=Z(\mathbb{F}_p,t)$ be the minimum number of field  operations required to determine all roots of any degree $t$ polynomial if it is known that all of them  are simple and lie in $\mathbb{F}_p$.  From \cite{BKS} we know that  $Z(t)\le tp^{1/2+o(1)}.$

Let  $R(F,s)$ be the minimum number of operations that are sufficient for the determination of the coefficients of the (minimum order)  recursion, of the given sequence $\{y_j\}_1^m\subset F$ if one knows that the sequence is given by a linear recursion of order not exceeding $s$ and $m>2s.$    The Berlekamp-Massey algorithm (see \cite{LN}, chapter  8)   after $O(s^2)$ field operations with $2s$ successive elements of a linear recurrent sequence of order  $s$ finds its characteristic polynomial.  Faster modifications of the algorithm are known: the one given in   (\cite{B}, chapter 11) requires   $O(s\log^{1+\eps} s) $ operations.

By  $LS(t)$ we denote the minimum number of operations in the field  $F$ sufficient to find solutions for any given non-homogeneous linear  $t\times t$ system of equations.  It is known that: $LS(t)=O(t^\omega),$ with the best current estimate  $\omega<2.4.$

\begin{theoremA}
There is a deterministic algorithm that takes as an input some vector  $y=\Phi x$ and  finds   after $O(m^{\omega}\log_p M+m^{2+o(1)}p^{1/2+o(1)})$  arithmetical operations in the field  $\mathbb{F}_p$ and      $O(m^2 \log_p M)$  operations in   $\mathbb{Z}$ the unknown vector $x$. Moreover the operations in $\Z$ are carried out with numbers not exceeding in absolute value  $\max\{M p,\max_j |y_j| +mp^{1-1/m}\}.$
\end{theoremA}

We describe  the strategy of the algorithm. We find step by step  the subsets  of the support $I_1\subseteq\dots I_\nu\dots \subseteq I$ of our vector  $x,$ and additionally we find finer  ``$p$--adic'' approximations for the coefficients $x_j$. Upper indices correspond to different steps of the algorithm.

1. Let $\alpha$ be the largest degree of  $p$, that divides all the coordinates of the vector  $y$, and $\beta$ the largest  degree of $p$, that divides all the coordinates of  $x.$ Let  $\tilde{y}=\frac{y}{p^{\alpha}}\in \Z^m, \tilde{x}=\frac{x}{p^{\beta}}\in \Z^p.$   It is clear that  (by the non degeneracy of the corresponding matrix over the field   $F$)  $\alpha=\beta.$ Therefore, the solution of our system is reduced to the solution of  $\Phi \tilde{x}=\tilde{y}.$  To simplify notation we suppose that $\alpha=\beta=0, x=\tilde{x},y=\tilde{y}.$

2. 
{\bf The main step of the algorithm.} We will analyze the system \eqref{sys1} in the field $F$. 
Let $p(t)=\prod_{j\in I} (t-j)=\sum_0^s p_l t^l\in F[t],$  where the coefficients $p_l$ are determined by  $I$ (and  $p_s=1$). It is easy to see that the sequence  $y_k$ is given by a linear recursion with characteristic polynomial $p$.  By a simple calculation  we have $\sum_{l=0}^s p_l y_{a+l}=0$  for all $a=0,1,\dots, m-s$

It is known that the minimum degree characteristic polynomial of the recursion divides the  characteristic polynomial of  any other recursion (\cite{LN}, chapter 8, theorem 8.42). Therefore its roots are in  $I$. We denote the set of roots by $S_1$.  Let  $I_1=S_1.$
 
 We start another iterative procedure  (a part of our main step). 
 
 {\bf Procedure.} In the field  $F$ we have the representation $y_l=\sum_{j\in I_1} \xi_j^{(1)}k_j j^l, 0\le l \le m-1.$  We find the characteristic polynomial of the minimum order recursion for  $\{y_l\}$, its roots and coefficients  $\xi_j^{(1)} \in F ,j\in I_1.$ 
 Let us note, that  $x_j=\xi_j^{(1)}, j\in I_1$ and $x_j=0, j\in I\setminus I_1$ (all equations are in $F$).  In fact, the vector from  $F^I$ with coordinates  $x_j-\xi_j^{(1)}, j\in I_1,$ and $x_j, j\in I\setminus I_1$  gives a solution to a non--degenerate homogeneous linear system by $\eqref{sys1}$  (its matrix is the Vandermonde matrix with columns multiplied by nonzero field elements).

Let us analyze the system of equations \eqref{sys1} in  $\mathbb{Z}.$ In what follows we identify different $\xi_j^{(\cdot)}\in F$ with appropriate  integers by taking the  minimum in absolute value representative of the residue class modulo $p$  i.e. from the set $\{0,\dots,\pm \frac{p-1}{2}\}.$  We compute the error vector  (with integer coordinates):    $ (y_l-\sum_{j\in I_1} \xi_j^{(1)}k_j j^l)_{0\le l \le m-1}$.  If this error vector  is zero, then  $I_1=I,$ the vector  $x$ is found and the algorithm stops.
Otherwise, there exists an integer  $\gamma_1\ge 1$  --- the degree of the highest power of  $p$ that divides every $y_l-\sum_{j\in I_1} \xi_j^{(1)} k_j j^l, 0\le l \le m-1$.   We define $x_j^{(1)}, j\in I$ and $y_l^{(1)},0\le l\le m-1$ as  $p^{\gamma_1}x_j^{(1)}=x_j$, $ j\in I\setminus I_1$,  $p^{\gamma_1}x_j^1=x_j-\xi_j^{(1)}$, $ j\in  I_1$ and $p^{\gamma_1} y_l^{(1)}=y_l-\sum_{j\in I_1} \xi_j^{(1)} k_j j^l, 0\le l \le m-1.$
From the considerations in 1 it follows that $x_j^{(1)},y_l^{(1)}\in \Z$ and the solution of  \eqref{sys1} reduces to the solution of   
\begin{equation}
\label{sys2}
\sum_{j\in I} k_j j^{l}x_j^{(1)}=y_{l}^{(1)},\,\,0\le l \le m-1.
\end{equation}
We try to find solutions in  $F$  for the system  $y_l^{(1)}=\sum_{j\in I_1} \xi_j^{(2)} k_j j^l, 0\le l \le m-1.$ While it is possible, we repeat the step of our procedure. In this way we make  $\kappa_1 $ steps of the procedure and obtain the system of equations  similar to  \eqref{sys2} equivalent to \eqref{sys1} in variables  $x_j^{(\kappa_1)}$.  But the  $\kappa_1+1$--step of the procedure is not possible.

 We apply the main step of the algorithm to the system of equations, obtained as a result of the procedure. We find the minimum degree polynomial of the recursion for $\{y_l^{(\kappa_1)}\}$; its roots make up $S_2\subset  I.$ It can happen that $I_1\cap S_2\not=\emptyset$.  But since we proceeded to the  $2$--nd step of the algorithm we find at least one new element of the support i.e. $S_2\setminus I_1\not = \emptyset.$ Let $I_2=I_1\cup S_2.$   After it we start the procedure once again...  The algorithm stops at some step (let it be   $L$).

3. Let us estimate the number of steps and arithmetical operations. The algorithm  has to stop since: at some step  $L$ for certain values of   $\kappa$ and $j$ we obtain  $x_j=t p^{\gamma_1+\dots +\gamma_\kappa}+\sum_{k=1}^{\kappa-1}\xi_j^{(k)}p^{\gamma_1+\dots +\gamma_k},$ where $t\in \Z\setminus\{0\}.$  It follows that  $|x_j|\ge \frac{1}{2}p^{\gamma_1+\dots +\gamma_\kappa}.$
We get $\gamma_1+\dots +\gamma_{\kappa_{L}}\le \log_p (2M),$ and therefore  $\kappa_L\le \log_p (2M).$
 Since, after each step of the algorithm we find at least one new element of the support  $I,$ the number of steps  $L\le s\le m.$ If the algorithm stops at the step number  $L$ then  the error vector equals  $0.$ We have $y_l^{(\kappa_L)}-\sum_{I_L} \xi_j^{(\kappa_L)} k_j j^l=0$ (in $\Z$) for all $0\le l\le m-1.$  There is also an element   $j\in I$ such, that  $\xi_j^{(\kappa_L)}\not=0.$  We remark that we have obtained the $p$--adic representation for  $x_j$. 
 
 We denote $s_\nu=|I_\nu|.$  On the step number  $\nu\in (1,L)$ we make no more than    $R(m)+Z(s_\nu)+(\kappa_\nu-\kappa_{\nu-1}) LS(s_\nu)$     arithmetical operations in  $F$ (determining  the characteristic polynomial for the recursion,  its zeroes, the coefficients  $\xi_j^{(\cdot)}$ (as a solution to the system of linear equations)),  and $ m(\gamma_{\kappa_{\nu-1}}+\dots +\gamma_{\kappa_{\nu}})+(\kappa_\nu-\kappa_{\nu-1})O(s_\nu)$ operations   (the determination of the quantity $\gamma_\nu$ (successive division on different powers of  $p$ the coordinates of the error vector), and the determination of  $y_l^{(\kappa_\nu)}$). We see that $\gamma_1+\dots+\gamma_{\kappa_L}\le \kappa_L\le \log_p (2M).$  Finally, the number of field operations is $O(m(m (\log m)^{1+o(1)} +mp^{1/2+o(1)})  +m^\omega\log_p M   )=O(m^{\omega}\log_p M+m^{2+o(1)}p^{1/2+o(1)})$ (we used the estimate $\omega> 2$), and number of the ring $\Z$ operations is  $O(m^2 \log_p M)$  (we include  the  $O(m\log_p M)$ operations for the determination of   $x$ using our $p$--adic approximations for its coordinates). Let us note that we make all computations in $\Z$ with numbers not exceeding    $\max\{M p,\max_j |y_j| +mp^{1-1/m}\}$ (in absolute value).

{\bf Remark.} There are different modifications of the algorithm. E.g. we can use other  algorithms for dealing with different problems in finite fields and the ring  $\Z$. If we apply the classical version of the Berlekamp-Massey algorithm (in order to find the recursion) the number of arithmetical operations in  $F$ can be estimated as    $O(m^{\omega}\log_p M+ m^3+m^2p^{1/2+o(1)}).$

 {\bf Acknowledgment.} I would like to express my deep gratitude to S.V.~Konyagin for his enthusiastic encouragement  and very useful discussions  concerning this research topic.


\begin{thebibliography}{99}

\bibitem{LN}
R. Lidl, H. Niederreiter {\it Finite fields.} Cambridge University  Press, 1984. 



\bibitem{FR}
S. Foucart, H. Rauhut,
{\it A Mathematical Introduction to Compressive Sensing}, Springer, New York, 2013.

\bibitem{FKS} L. Fukshansky, D. Needell, and B. Sudakov,
An algebraic perspective on integer sparse recovery,
{\it Applied Mathematics and Computation} {\bf 340} (2019), 31--42.

\bibitem{K}  S. V. Konyagin, On the recovery of an integer
vector from linear measurements, {\it Mathematical Notes},
{\bf 104} (2018), 859--865.

\bibitem{KS}
S.V. Konyagin, B. Sudakov,  An extremal problem for integer sparse recovery, arXiv: 1904.08661

\bibitem{BKS} J. Bourgain, S.V. Konyagin, I.E. Shparlinski,  Character sum bounds and deterministic polynomial root finding in finite fields, {\it Mathematics of Computation}, {\bf 84(296)} (2015), 2969-2977.

\bibitem{B}
R.E. Blahut, {\it Theory and Practice of Error Control Codes}, Addison- Wesley, 1983.



\end{thebibliography}
\end{document}